\documentclass{IEEEtran} 
\usepackage{cite}
\usepackage{amsthm}
\usepackage{amsmath,amssymb}
\usepackage{subcaption}
\usepackage{float}
\usepackage{bbm}
\usepackage{listings}
\usepackage{graphicx}
\usepackage{amsfonts}
\usepackage{color}
\usepackage{enumerate}
\usepackage{bm}
\usepackage{stmaryrd}
\usepackage{comment}
\usepackage{cancel}
\usepackage{soul}
\usepackage{tikz}
\usetikzlibrary{angles,quotes}

\newtheorem{theorem}{Theorem}

\newtheorem{proposition}{Proposition}

\newtheorem{definition}{Definition}

\def\BibTeX{{\rm B\kern-.05em{\sc i\kern-.025em b}\kern-.08em
    T\kern-.1667em\lower.7ex\hbox{E}\kern-.125emX}}

\graphicspath{{./}}

\input{CS_Hilbert.sty}

\input{tcilatex}

\begin{document}
\title{Reconstructing high-dimensional Hilbert-valued functions via compressed sensing}
\author{Nick Dexter, Hoang Tran, and Clayton Webster 
\thanks{This material is based upon work supported in part by: the U.S. Department of Energy, Office of Science, Office of Advanced Scientific Computing Research, Applied Mathematics program under contracts and awards ERKJ314, ERKJ331, ERKJ345, and Scientific Discovery through Advanced Computing (SciDAC) program through the FASTMath Institute under Contract No. DE-AC02-05CH11231; and by the Laboratory Directed Research and Development program at the Oak Ridge National Laboratory, which is operated by UT-Battelle, LLC., for the U.S. Department of Energy under contract DE-AC05-00OR22725.}
\thanks{Department of Mathematics, Simon Fraser University. Email: nicholas$\_$dexter@sfu.ca}
\thanks{Computer Science and Mathematics Division, Oak Ridge National Laboratory. Email: tranha@ornl.gov }
\thanks{Department of Mathematics, University of Tennessee, Knoxville, and Computer Science and Mathematics Division, Oak Ridge National Laboratory. Email: webstercg@math.utk.edu}}

\maketitle

\begin{abstract}
We present and analyze a novel sparse polynomial technique for approximating high-dimensional Hilbert-valued functions, with application to parameterized partial differential equations (PDEs) with deterministic and stochastic inputs. Our theoretical framework treats the function approximation problem as a joint sparse recovery problem, where the set of jointly sparse vectors is possibly infinite. To achieve the simultaneous reconstruction of Hilbert-valued functions in both parametric domain and Hilbert space, we propose a novel mixed-norm based $\ell_1$ regularization method that exploits both energy and sparsity. Our approach requires extensions of concepts such as the restricted isometry and null space properties, allowing us to prove recovery guarantees for sparse Hilbert-valued function reconstruction.
We complement the enclosed theory with an {algorithm for Hilbert-valued recovery, based on standard forward-backward algorithm}, meanwhile establishing its strong convergence in the considered infinite-dimensional setting. 
Finally, we demonstrate the minimal sample complexity requirements of our approach, relative to other popular methods, with numerical experiments approximating the solutions of high-dimensional parameterized elliptic PDEs.
\end{abstract}


\section{Introduction}

\label{sec:intro}

Underlying many successful applications of compressed sensing to problems in applied mathematics and the physical sciences is the fact that, for many practical problems, the target to be reconstructed possesses sufficient sparsity in order to enable unique solutions from systems that would otherwise be ill-posed.
In the basic compressed sensing problem, the signal is an unknown vector $\bc \in \mathbb{R}^{N}$, and the sensing process yields a measurement vector $\bu \in \mathbb{R}^m$ that is formed by the product of $\bc$ with a sensing matrix, i.e., $\bu = \bA \bc$, where $\bA \in \mathbb{R}^{m\times N}$. 
The key observation is that when the signal $\bc$ is sufficiently sparse, it can be uniquely determined from an underdetermined set of measurements ($m\ll N$), provided $\bA$ satisfies certain additional properties.
To overcome the NP-hardness of directly finding the sparsest $\bc$ consistent with a given set of measurements, various greedy and convex relaxation strategies have been proposed and demonstrated, both empirically and theoretically, to have good reconstruction performance in a range of settings. 

In the parameterized PDE literature, one typically seeks to approximate the parameter-to-solution map $\by \mapsto u(\by)$, taking values in a Hilbert space $\cV$, by a truncation of its orthonormal expansion, i.e., 
\begin{align}
\label{eq:orthonormal_truncation}
u = \sum_{\bmnu\in\cF} \bc_{\bmnu} \Psi_\bmnu \approx \sum_{\bmnu \in \cJ} \bc_\bmnu \Psi_\bmnu =: u_{\cJ}.
\end{align}
Here, $\by\in \mathcal{U}$ a high-dimensional tensor product domain, $\cJ$ is a finite multi-index set of cardinality $N=\#(\cJ)$ with $\cJ \subset \cF := \N_0^d$, $(\Psi_\bmnu)_{\bmnu\in\mathcal{J}}$ is an {$L^2(\mathcal{U},d{\varrho})$}-orthonormal basis, and $\bc = (\bc_\bmnu)_{\bmnu\in\cJ} \in \mathcal{V}^N$ are the Hilbert-valued coefficients to be computed.
Generally, one selects $\cJ$ to be a set large enough to ensure that 
\begin{align}
\label{eq:truncation_error}
\|u - u_{\cJ} \|_{L^2}^2  = \| u_{\cJ^c}\|_{L^2}^2 = \sum_{\bmnu\in\cF\setminus\cJ} \| \bc_{\bmnu} \|_{\cV}^2
\end{align}
is minimal. However, this can lead to a less efficient approximation if $\cJ$ is not chosen carefully.

Over the course of the last decade, a series of works in the parameterized PDEs community (see \cite{CohenDeVore15} and the reference therein)
have demonstrated that, under reasonable assumptions on the input data to the PDEs, the solutions are compressible, and hence well-represented by sparse expansions in given orthonormal systems.
In other words, the solution vector $\bc$ from \eqref{eq:orthonormal_truncation} is sparse, and accurate reconstructions of the most important components of $\bc$ are enough for satisfactory approximations of the PDE solutions. The ability of compressed sensing (CS) to exploit sparsity and allow far fewer samples than traditional approaches (e.g., Monte-Carlo, projection, interpolation) makes it a promising tool for such reconstruction problems. Therefore it is no surprise that 
CS-based polynomial approximation has attracted growing interest in the area of high-dimensional parameterized PDEs in recent years, \cite{DO11, MG12,  YK13, Rauhut2017, PHD14, HD15, PengHampDoos15}. 

However, there has been a critical mismatch between standard CS techniques and the problem of reconstructing Hilbert-valued
solutions to parameterized PDEs: these methods do not enable direct recovery of the vector $(\bc_\bmnu)_{\bmnu\in\cJ}\in\cV^N$, i.e., a vector with Hilbert-valued components. 
Instead, CS-based polynomial approximation methods only allow the recovery of real or complex sparse vectors. 
Hence, 
in the context of parameterized PDEs,
standard CS-based polynomial approximation methods do not yield approximations of the entire solution map $\by \mapsto u(\cdot, \by) \in \cV$, but only functionals of the solution, i.e., maps of the form $\by \mapsto G(u(\by)) \in \mathbb{R}$, where $G: \cV \to \mathbb{R}$ is a functional of $u$. In particular, many of the existing works (cited above) perform \textit{point-wise reconstruction} of the solution, i.e., reconstruction of the map $\by \mapsto G_{x^\star}(u(\by))$ with $G_{x^\star}(u(\by)) = u(x^{\star},\by) $ for a fixed $x^{\star} \in D$, with $D$ the physical domain of the PDE. 

In this paper, we {present an overview of latest progress in addressing} this mismatch by developing a novel sparse recovery technique which enables direct reconstruction of Hilbert-valued vectors, {\cite{DexterTranWebster2018,DexterTranWebster17}}. The theoretical framework, based on the problem of joint-sparse recovery, treats the recovery problem as a matrix recovery problem in which each row of the solution matrix may have infinitely many terms corresponding to the coefficients of a Hilbert-valued function in a given basis.
Our regularization is performed with respect to a mixed norm {$\|\cdot\|_{\cV,1}$}, which is defined to be the {$\ell_1$} norm of the vector $(\|\bc_\bmnu\|_{\cV})_{\bmnu\in\cJ}$. Often, the decay of the polynomial coefficients and tail {bounds} are estimated in global energy norms rather than pointwise over the physical domain $D$ \cite{CohenDeVore15}. As we will show, our choice of regularization enables us to prove convergence rates for global approximations in {the same energy} norms, with explicit computable coefficients. We present the theoretical support of this strategy via several extensions of compressed sensing concepts such as the {restricted isometry property} (RIP) and {null space property} (NSP) to the Hilbert-valued setting. 

Our regularization problem is solved with a forward-backward splitting approach. We present a new strong convergence result for forward-backward splitting in a joint sparse recovery scenario where we assume neither the strict convexity of the {fidelity} functions nor the finite cardinality of the set of signals to construct. We note that most of similar strong convergence results along the line need at least one of these conditions. In relation to parameterized PDEs problem, this result means that forward-backward splitting strongly converges even before discretization in (infinite-dimensional) Hilbert space is introduced.

The rest of the paper proceeds as follows. In Section \ref{sec:parameterized_PDE}, we introduce our approach of sparse regularization and provide a brief comparison to the problem of joint-sparse recovery. 
In Section \ref{sec:recovery}, we present theoretical guarantees for the direct reconstruction of Hilbert-valued functions through mixed norm regularization, as well as convergence estimates in case of parameterized PDEs. Section \ref{sec:forward_backward} introduces a version of forward-backward algorithm adapted for Hilbert-valued recovery and provides its strong convergence result. 
Section \ref{sec:numerical} presents numerical results on applying the sparse regularization method to the solution of a parameterized elliptic PDE.
Finally, Section \ref{sec:conclusion} concludes with some remarks on our approach and applications to a wider array of problems.

\section{Sparse regularization for parameterized PDEs}

\label{sec:parameterized_PDE}

The sparse polynomial techniques proposed in this work are applicable to general parameterized PDE problems of the form: find $u(\cdot,\bm{y}) : \overline{D} \times \mathcal{U} \to \mathbb{R}$ for all $\by \in \mathcal{U}$ such that  
\begin{equation}
\label{eq:genPDE}
\mathcal{L}(u(\cdot, \bm{y}),\bm{y})=0, \ \ \text{ in }D, 
\end{equation}
where $\mathcal{L}$ is a differential operator defined on a spatio-temporal domain $D$.
Our development to the general CS/joint-sparse polynomial approximation problem can be formulated as follows.
One first generates $m$ samples $\by_1, \ldots, \by_m$ in $\cU$ independently from the orthogonalization measure $\varrho$ associated with $(\Psi_\bmnu)_{\bmnu \in \cJ}$, {for instance, uniform samples for the Legendre basis and Chebyshev samples for Chebyshev basis, see \cite[Chapter 12]{FouRau13},} and solves the equation \eqref{eq:genPDE} at these samples to form the normalized output vector $\bu := (u(\by_1),\ldots,u(\by_m))/\sqrt{m}$ as well as the normalized sampling matrix $\bm{A} := (\Psi_{\nu}(\by_i)/\sqrt{m})_{{{1\le i \le m},\, {\nu \in \cJ}}}$. 
Taking account that the true unknown coefficient $\bc = (\bc_{\bmnu})_{\bmnu \in \cJ}$ approximately solves the linear system 
\begin{align}
\label{und_system}
\bu = \bm{A} \bz, \ \ \ \bz\in \cV^N,
\end{align}
and further, $\bc$ is compressible \cite{CohenDeVore15}, it is reasonable to approximate $\bc$ by $\bc^\#$, the solution to the following problem
\begin{align}
\label{eq:lV1_BPDN}
\textnormal{min}_{\bz\in\cV^N} \; \|\bz\|_{\cV,1} \;\; \textnormal{subject to} \;\; \|\bA \bz - \bu\|_{\cV,2} \le \frac{\eta}{\sqrt{m}},
\end{align}
{where $\eta$ relates to an estimate of the tail \eqref{eq:truncation_error}}; or the equivalent unconstrained convex minimization:
\begin{align}
\label{eq:lV1}
\textnormal{min}_{\bz \in {\cV}^N} \; \|\bz\|_{\cV,1} + \frac{\mu}{2}  \|\bm{A} \bz - \bu\|^2_{\cV,2},
\end{align}
for appropriately chosen $\mu>0$.
Here, the norm $\|\cdot\|_{\cV,q}$ is defined for $\bc\in \cV^N$ as $\|{\bc}\|_{\cV,q} := (\sum_{\bmnu \in \cJ} \|\bc_{\bmnu}\|^q_{\cV})^{1/q}$. This is arguably the most natural extension of the $\ell_1$ minimization approach, traditionally for real and complex signal recovery, to the reconstruction of sparse generalized Hilbert-valued vectors. We denote our approach {\em simultaneous compressed sensing} (SCS). In depth description, analysis and application of SCS for solving parameterized PDEs are provided in \cite{DexterTranWebster17}.  

Problem \eqref{eq:lV1_BPDN} can be related to the joint-sparse basis pursuit denoising problem as follows.
Let $(\phi_r)_{r\in \mathbb{N}}$ be an orthonormal basis of $\cV$, then $\bc_{\bmnu} \in \cV$ has unique representation 
$$
\bc_\bmnu  = \sum_{r\in \mathbb{N}} c_{\bmnu,r} \phi_r,\ \ \text{with } \ c_{\bmnu,r} \in \mathbb{R}.  
$$
Each coefficient $\bc_{\bmnu}$ corresponds to an $\mathbb{R}^{1\times \mathbb{N}}$ vector $(c_{\bmnu,1},c_{\bmnu,2},\ldots, c_{\bmnu,r}, \ldots)^{\top}$, thus, $\bc = (\bc_{\bmnu})_{\bmnu \in \cJ}$ is completely determined by the $\mathbb{R}^{N\times \mathbb{N}}$ matrix ${\widehat{\bm{c}}} = (c_{\bmnu,r})_{\bmnu \in \cJ,r\in \mathbb{N}}$. Furthermore, 
\begin{align*}
\|\bc\| _{\cV,q} = \left(\sum_{\bmnu \in \cJ} \|\bc_{\bmnu}\|_{\cV}^q \right)^{1/q} & = \left(\sum_{\bmnu \in \cJ} \left( \sum_{r \in \mathbb{N}} |c_{\bmnu,r}|^2\right) ^{q/2}\right)^{1/q} \\
& \equiv  \left(\sum_{\bmnu \in \cJ} \|\widehat{\bc}_{\bmnu}\|_2^q\right)^{1/q}  = \|\widehat{\bc}\|_{2,q},
\end{align*}
where the matrix norm $\|\cdot\|_{p,q}$ is defined as the $\ell_{q}$ norm of the vector $(\|\widehat{\bc}_{\bmnu}\|_{p})_{\bmnu\in\cJ}$, implying the equivalence of \eqref{eq:lV1_BPDN} with the infinite-dimensional joint-sparse recovery problem.
In practice, one needs to employ a discretization over $\mathcal{V}$ to be able to numerically solve \eqref{eq:lV1_BPDN} or \eqref{eq:lV1}. Any preferred method of spatial/temporal discretization may be used, e.g., the finite element, difference, and volume methods. However, our strong convergence result of forward-backward splitting in Section \ref{sec:forward_backward} is applicable to \eqref{eq:lV1} without any discretization.

\section{Error estimates for Hilbert-valued recovery}
\label{sec:recovery}

Straightforward extensions of concepts and results from compressed sensing and joint-sparse recovery can be made to ensure uniform recovery of Hilbert-valued signals via $\ell_{\cV,1}$-relaxation.
Well-known concepts such as the {null space property} (NSP) and {restricted isometry property} (RIP) have Hilbert-valued counterparts.
In this section, we review Hilbert-valued versions of the NSP and RIP to guarantee uniform recovery for the mixed norm regularization in the presence of noise or sparsity defects. 
{We note that the extension of these concepts to the Hilbert-valued recovery setting does not require finite-dimensionality of the Hilbert space $\cV$. For more detailed discussion and complete proofs of the results in this section, we refer the interested readers to \cite{DexterTranWebster17}.}

{First, we define a Hilbert-valued version of the $\ell_2$-robust NSP, which guarantees the reconstruction of vectors ${\bm c} \in \mathcal{V}^N$ (up to $s$ largest components and up to a noise level).}

\begin{definition}[$\ell_{\cV,2}$-robust null space property]
\label{def:l_V2_RNSP}
The matrix $\bA\in\R^{m\times N}$ is said to satisfy the $\ell_{\cV,2}$-robust null space property of order $s$ with constants $0< \rho < 1$ and $\tau > 0$ if 
\begin{align}
&\|\bz_S \|_{\cV,2} \le \frac{\rho}{\sqrt{s}} \|\bz_{S^c}\|_{\cV,1} + \tau \|\bA \bz \|_{\cV,2}
\\
 &\qquad\qquad\quad\quad \forall \bz\in\cV^N, \forall S\subset[N] \;\; \textnormal{with} \;\; \#(S) \le s. \notag
\end{align}
\end{definition}

An RIP-type condition is required to quantify the sample complexity of solving \eqref{eq:lV1_BPDN} to a prescribed accuracy.
The following result establishes the implication of the $\ell_{\cV,2}$-robust NSP from the standard RIP.

\begin{proposition}
\label{prop:RIP_V_RIP_equivalence}
Suppose that $\cV$ is a separable Hilbert space, and that the matrix $\bA\in \R^{m\times N}$ satisfies the RIP, that is
\begin{align}
\label{eq:standard_RIP}
&(1 - \delta_{2s}) \|\bz \|_2^2 \le \|\bA \bz \|_2^2 \le (1+ \delta_{2s}) \|\bz \|_2^2, \\
&\qquad\qquad\quad\quad \forall \bz\in\R^N, \;\; \bz \;\; 2s\textnormal{-sparse}, \notag
\end{align}
with $\delta_{2s} < \frac{4}{\sqrt{41}}$. Then, $\bA$ satisfies the $\ell_{\cV,2}$-robust NSP of order $s$ with constants $0< \rho < 1$ and $\tau>0$ depending only on $\delta_{2s}$.
\end{proposition}

Proposition \ref{prop:RIP_V_RIP_equivalence} implies that sample complexity results for solving the standard $\ell_1$ problem
hold for the mixed norm problem \eqref{eq:lV1_BPDN} as well. 

In the parameterized PDE context, with error estimates as in \cite[Theorem 2]{TWZ17} for quasi-optimal approximations, we are able to provide convergence rates for approximations to parameterized PDEs obtained through the mixed norm regularization. 
First, we assume that the solution $u$ has parametric expansion with coefficients $(\bc_{\bmnu})_{\bmnu\in\cF}$ as in \eqref{eq:orthonormal_truncation} satisfying $\|\bc_{\bmnu}\|_\cV \lesssim e^{-b(\bmnu)}$ for every $\bmnu\in\cF$, with $b(\bmnu)$ obeying \cite[Assumption 3]{TWZ17}. For brevity, we do not detail that assumption herein, but remark that it is satisfied by most parameterized PDE models we are aware of.

\begin{theorem}
\label{thm:SCS_convergence}
For any $\varepsilon > 0$, assume that the solution $u$ to \eqref{eq:genPDE} with parametric expansion \eqref{eq:orthonormal_truncation} has coefficients $(\bc_{\bmnu})_{\bmnu\in\cF}$ satisfying $\|\bc_\bmnu\|_\cV \le e^{-b(\bmnu)}$ for every $\bmnu\in\cF$ with $b:[0,\infty)^d \to \R$ also satisfying \cite[Assumption 3]{TWZ17}.
Denote by $\cJ_s$ the set of indices corresponding to the $s$ largest bounds of the sequence $(e^{-b(\bmnu)})_{\bmnu\in\cF}$.
Let $\cJ$ be such that $\cJ_s \subseteq \cJ$, and assume that the number of samples $m\in\N$ satisfies
\begin{align}
\label{eq:RIP_complexity}
&m\ge C \Theta^2 s \max \{\log^2(\Theta^2 s) \log(N),
\\
&\qquad\qquad\qquad\qquad \log(\Theta^2 s) \log(\log(\Theta^2 s) N^{\log(s)}) \}, \notag
\end{align}
with $\Theta = \sup_{\bmnu\in\cJ} \|\Psi_\bmnu\|_{L^\infty(\cU)}$ and $N= \#(\cJ)$. 
Then with probability $1-N^{-\log(s)}$, the solution $\bc^\#$ of 
\begin{align}
    \label{eq:SCS_BPDN_repeat}
\textnormal{minimize}_{\bz\in\mathcal{V}^N} \;\; \|\bz\|_{\mathcal{V},1} \quad \textnormal{s.t.} \quad \|\bA \bz - \bm{u} \|_{\mcV,2} \le \frac{\eta}{\sqrt{m}}
\end{align}
approximates $u$ with asymptotic error
\begin{align}
\label{eq:SCS_convergence}
\left\| u - \sum_{\bmnu \in \cJ} \bc_\bmnu^\# \Psi_\bmnu \right\|_{L^2_\varrho(\cU;\cV)} \le {C}_\varepsilon \sqrt{s} \exp \left[ - \left( \frac{\kappa \, s}{(1 + \varepsilon)} \right)^{1/d} \right],
\end{align}
where $\kappa,{C}_\varepsilon >0$ are independent of $s$.
\end{theorem}

\section{Forward-backward algorithm for Hilbert-valued recovery}

\label{sec:forward_backward}

In this section, we present a forward-backward splitting algorithm for solving problem \eqref{eq:lV1} over the real Hilbert space $\cV^N$. 
Assuming $\mu = 1$ for simplicity, define $T_1$ and $T_2$ to be the subdifferential and gradient of the $\|\cdot\|_{\cV,1}$ and $\frac{1}{2} \|\bA(\cdot) - \bu\|_{\cV,2}^2$ parts of \eqref{eq:lV1}, respectively, and $T=T_1+T_2$.
The algorithm can be derived as follows. Given $\tau > 0$, we have
\begin{align*}
\bm{0} \in T(\bx) 
                  & \iff \bx = (I + \tau T_1)^{-1} (I - \tau T_2) \bx. \numberthis \label{eq:forward_backward_derivation}
\end{align*}
The last identity in \eqref{eq:forward_backward_derivation} leads to the {\em forward-backward splitting} algorithm: 
given initial guess $\bx^{0}\in \mathcal{H} := \{\bx \in \cV^{N} : \|\bx\|_{\cV,2} < \infty\}$, compute
\begin{equation}
\label{eq:forward_backward_iteration}
\bx^{k+1} = (I + \tau T_1)^{-1} (I - \tau T_2)\, \bx^k, 
\end{equation}
where $\bx^{k}$ denotes the approximation at $k$-th iterate. 
Letting
$$
J_\tau := (I + \tau T_1)^{-1},\ \ G_\tau := (I - \tau T_2),\ \ \text{ and }\ S_\tau := J_\tau \circ G_\tau,
$$
then \eqref{eq:forward_backward_iteration} can be written as $\bx^{k+1}  = S_\tau (\bx^{k})$.
A straightforward derivation yields: 
\begin{align}
\label{eq:G_tau}
G_\tau(\bx) & = \bx - \tau \bA^*(\bA\bx - \bu), \\
\label{eq:J_tau}
J_{\tau,j}(\bx) 
 & = \frac{\bx_{j}}{\|\bx_{j}\|_2} \cdot \max \{ \|\bx_{j} \|_2 - \tau, 0\}, \quad 1\le j\le N. 
\end{align}
One can observe that the forward operator $G_\tau$ resembles a step of gradient descent algorithm with stepsize $\tau$. 
The backward operator $J_{\tau}$, on the other hand, is a soft  thresholding step associated with proximal point method. 

Under a standard assumption on the stepsize $\tau$, namely $0<\tau < 2/\|\bA^*\bA\|_{2}$, we can obtain some nonexpansive properties for $G_\tau$ and $J_\tau$. 
In particular, $G_\tau$ is nonexpansive, i.e., 
\begin{align}
\label{eq:G_tau_FNE}
\| G_\tau (\bv) - G_\tau (\bw) \|_{2,2} \le \| \bv - \bw \|_{2,2} \qquad \forall \bv,\bw\in \cH,
\end{align}
and $J_\tau$ is row-wise firmly nonexpansive, i.e., $\forall \bv, \bw \in \cH$, $\forall j \in [N] := \{1,2,\ldots,N\}$ 
  \begin{gather}
  \label{eq:J_tauj_FNE}
\begin{aligned}
 & \| J_{\tau} (\bv_j) - J_{\tau} (\bw_j) \|_{2}^2 \\
 & \;\; \le \| \bv_j - \bw_j \|_{2}^2 - \| (I - J_{\tau}) \bv_j  - (I - J_{\tau}) \bw_j \|_{2}^2. 
\end{aligned}
\end{gather}
see \cite[Chapter 4]{Bauschke2010}. {To prove the strong convergence of the forward-backward splitting algorithm for Hilbert-valued recovery in infinite-dimensional setting, we consider certain partition of the index set $[N]$ into two subsets $L$ and $E$, $L\subset (\text{supp}(\bm x^*))^c$ and $E\supset \text{supp}(\bm x^*)$ where $\bx^*$ is a solution to \eqref{eq:lV1}, as they require different treatments, see \cite{HYZ08,DexterTranWebster2018}.  
With properties \eqref{eq:G_tau_FNE}-\eqref{eq:J_tauj_FNE} in hand, our result is obtained in three steps:}

\begin{enumerate}

\item Using a partitioning technique from \cite{HYZ08}, establish finite convergence $\bx_j^k \to \bx^*_j$ for every $j\in L$.

\item Establish an angular convergence result on the extended support set $E$. 

\item Combine the known general weak convergence of the forward-backward algorithm, see, e.g., \cite{Bauschke2010}, with our angular convergence result to obtain convergence in norm on $E$, and hence strong convergence since $L\cup E = [N]$.
\end{enumerate}

For more details on the above approach, see \cite{DexterTranWebster2018}. The established result can then be summarized as follows.

\begin{theorem}
\label{thm:strong_convergence}
Let $0 < \tau < 2/\|\bA^*\bA\|_{2}$ and $(\bx^{k})$ be the sequence generated by the forward-backward iterations $\bx^{k+1} = S_\tau (\bx^k)$ starting from any $\bx^{0}\in \cH$. Then $(\bx^{k})$ converges strongly to an element $\bx^*\in \cV^N$ solving \eqref{eq:lV1}.
\end{theorem}

\section{Numerical experiments on parameterized elliptic PDE models}

\label{sec:numerical}

In this section we present numerical experiments demonstrating the efficiency of the proposed approach in approximating the solution of the following parameterized elliptic PDE: for all $\by\in\cU$, find $u(\cdot,\by):D \to \R$ such that
\begin{align}
\label{eq:model_problem}
\left\{\begin{array}{rll} 
-\nabla \cdot \left( a(x,\by) \nabla u(x,\by) \right) \hspace{-0.25cm} &= f(x)  & \textnormal{in } D, \\
                                                u(x,\by) \hspace{-0.25cm} &= 0  & \textnormal{on } \partial D.
\end{array}\right.
\end{align}
Here $D=[0,1]^2$,  $f \equiv 1$, the diffusion coefficient is given by
\begin{align*}
a(x,\by) & = 10 + y_1 \left(\frac{\sqrt{\pi}L}{2}\right)^{1/2} + \sum_{i=2}^d \; \zeta_i \; \vartheta_i(x) \; y_i \numberthis \label{eq:transcendental_a_linear_version} \\
\zeta_i & := (\sqrt{\pi} L)^{1/2} \exp\left( \frac{-\left( \left\lfloor \frac{i}{2} \right\rfloor \pi L\right)^2}{8} \right), \text{ for } i>1,  \\
\vartheta_i(x) & := \left\{ \begin{array}{rl} \sin\left( \left\lfloor\frac{i}{2}\right\rfloor \pi x_1/L_p \right),  \text{ if $i$ is even,} \\ \cos\left(\left\lfloor\frac{i}{2}\right\rfloor \pi x_1/L_p \right),  \text{ if $i$ is odd.} \end{array} \right.
\end{align*}
or its log transformed version $\log(a(x,\by) - 0.5)$, and $y_i$ are uniform random variables on $(-\sqrt{3},\sqrt{3})$.
Convergence is compared against ``highly-enriched'' reference sparse grid stochastic collocation approximations which we denote $u_{h,\textnormal{ex}}(x,\by)$, see \cite{Stoyanov:2016ci}. 
All approximations (including the enriched reference solution) are computed on fixed finite element meshes $\mathcal{T}_h$, and our enriched SC approximation is computed using {\em Clenshaw-Curtis} (CC) abscissas with high level $L_{\textnormal{ex}}$.
We compare performance using the relative errors of the expectation and standard deviation in the $H_0^1(D)$-norm. 

In our plots and discussion, we use the following abbreviations. For the SCS method, we use ``SCS-TD'' to denote approximations obtained using the SCS method with the {\em total degree} (TD) subspace. For the {\em stochastic Galerkin} (SG) method, we use ``SG-TD'' to denote the SG approximation with TD subspace, see, e.g., \cite{Xiu2002,Gunzburger:2014hi}. The {\em stochastic collocation} (SC) method with CC points and doubling growth rule is denoted ``SC-CC,'' see \cite{NTW08}. The Monte Carlo approximation is denoted ``MC.'' 

We follow convention from \cite{BNTT_comp}, identifying {\em stochastic degrees of freedom} (SDOF) as the number of random sample points $m$ for MC and SCS, and sparse grid points $m_L$ for SC with level $L$. 
For SG, we define SDOF to be $N$, the cardinality of the index set used in construction. 
We include the SG method only to compare the $L^2(\cU,d\varrho)$-optimal (w.r.t. SG SDOF) error of the Galerkin projection against the error of the sampling-based approximations. 

We employ the orthonormal Legendre series for the parametric discretization of the SG and SCS methods. 
For the MC and SCS methods, the random samples $(\by_i)_{i=1}^m$ are drawn from the uniform distribution $\varrho$. 
We average the random sampling results over 24 trials, fixing the initial seed for the pseudorandom number generator on each trial, and then solving each trial's problem with the same set of $m$ samples when plotting convergence.

\begin{figure}[ht]
\begin{center}
\includegraphics[clip=true,trim=00mm 00mm 10mm 00mm,width=0.24\textwidth]{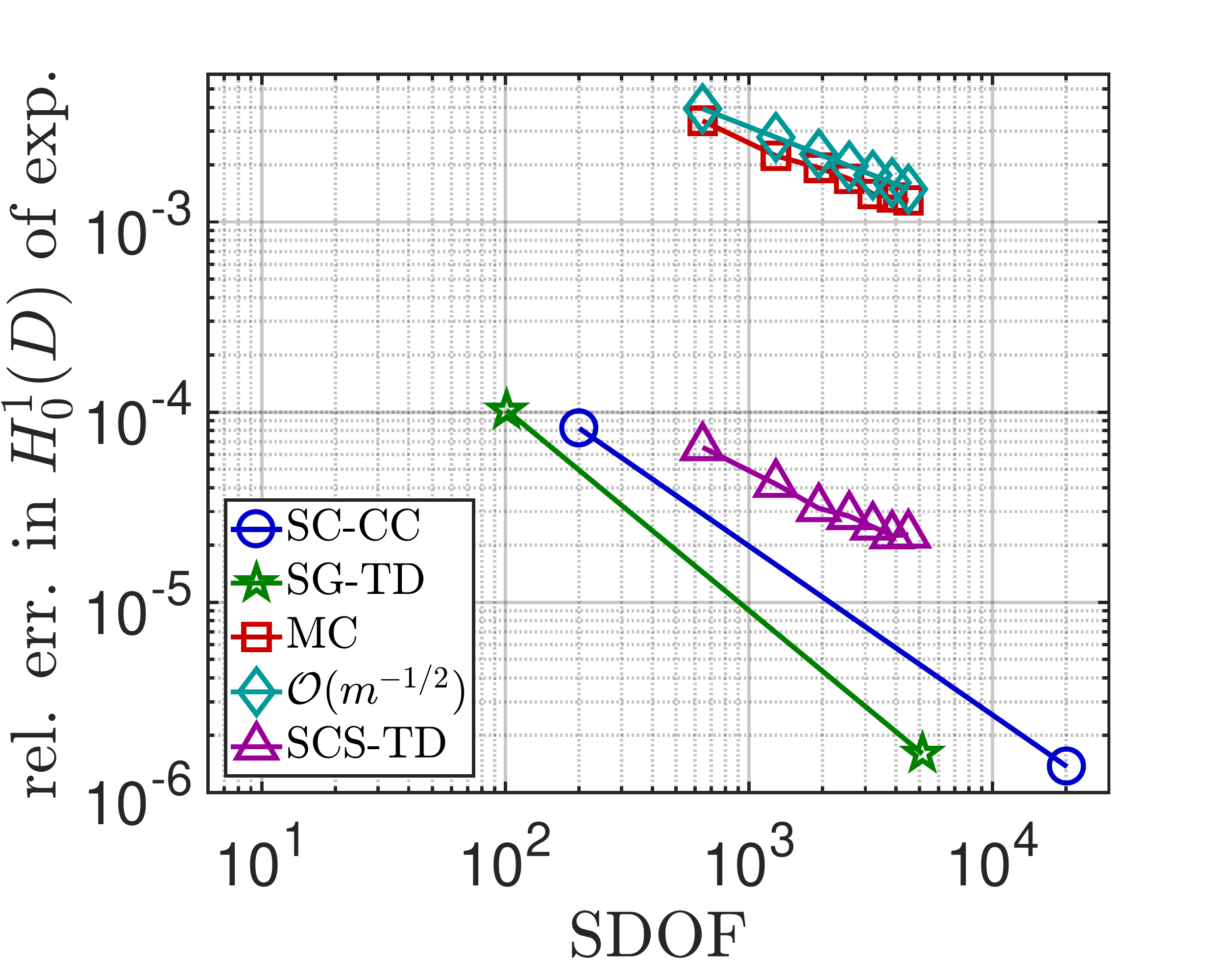}
\includegraphics[clip=true,trim=00mm 00mm 10mm 00mm,width=0.24\textwidth]{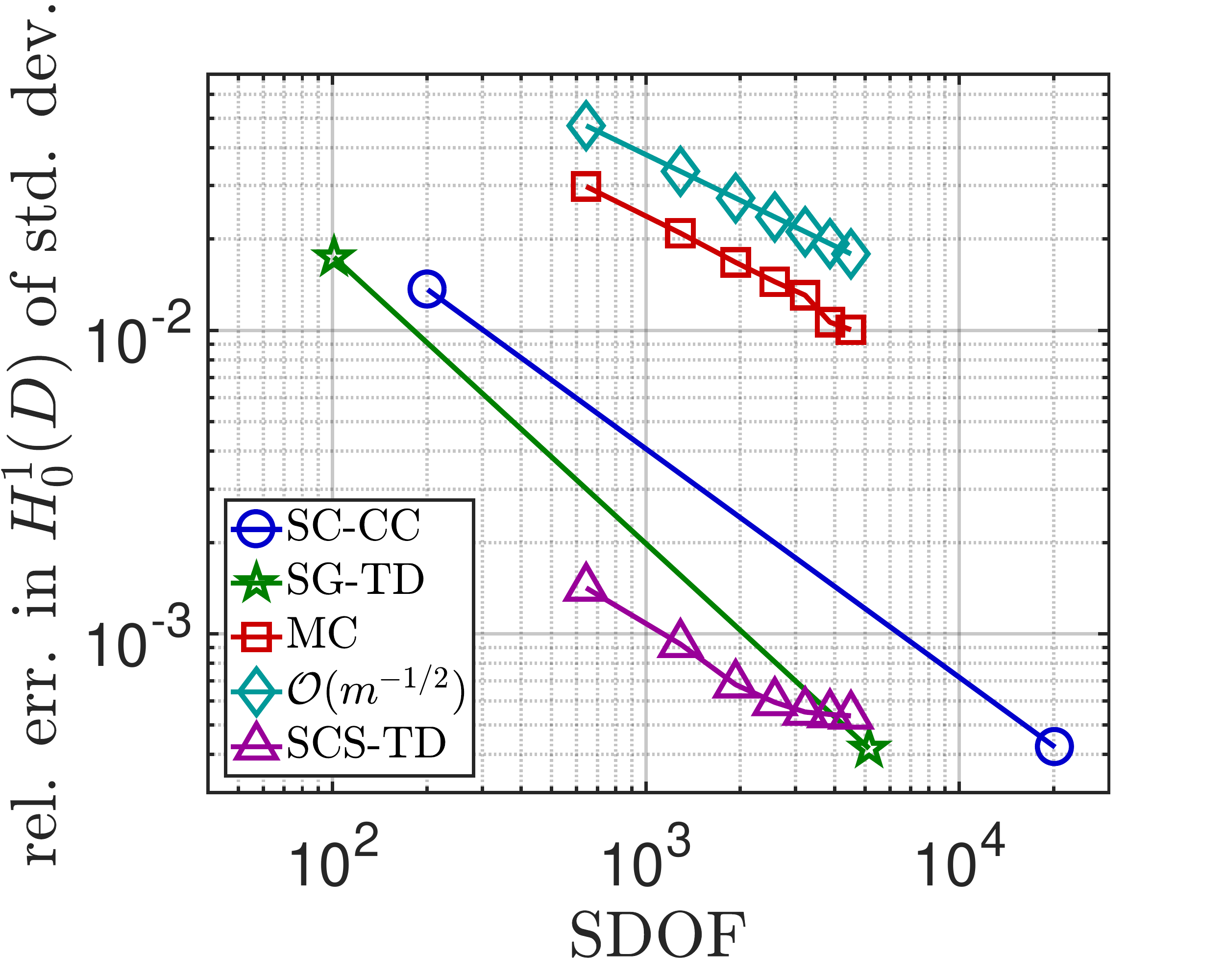}
\end{center}
\caption{Comparison of relative errors in expectation \textbf{(left)} and standard deviation \textbf{(right)} in the $H_0^1(D)$ norm for the $L=2,3$ stochastic collocation (SC-CC),  $p=1,2$ stochastic Galerkin (SG-TD), Monte Carlo (MC), and total degree order $p=2$ with $N=5151$ simultaneous compressed sensing (SCS-TD) methods for solving \eqref{eq:model_problem} with coefficient \eqref{eq:transcendental_a_linear_version} and correlation length $L_c=1/2$ in $d= 100$ dimensions. All approximations computed on a finite element mesh with 206 degrees of freedom corresponding to a maximum mesh size of $h\approx 1/16$.}
\label{fig:CS_vs_others_increasing_correlation_length}
\end{figure}

\begin{figure}[ht]
\begin{center}
\includegraphics[clip=true,trim=00mm 04mm 14mm 00mm,width=0.24\textwidth]{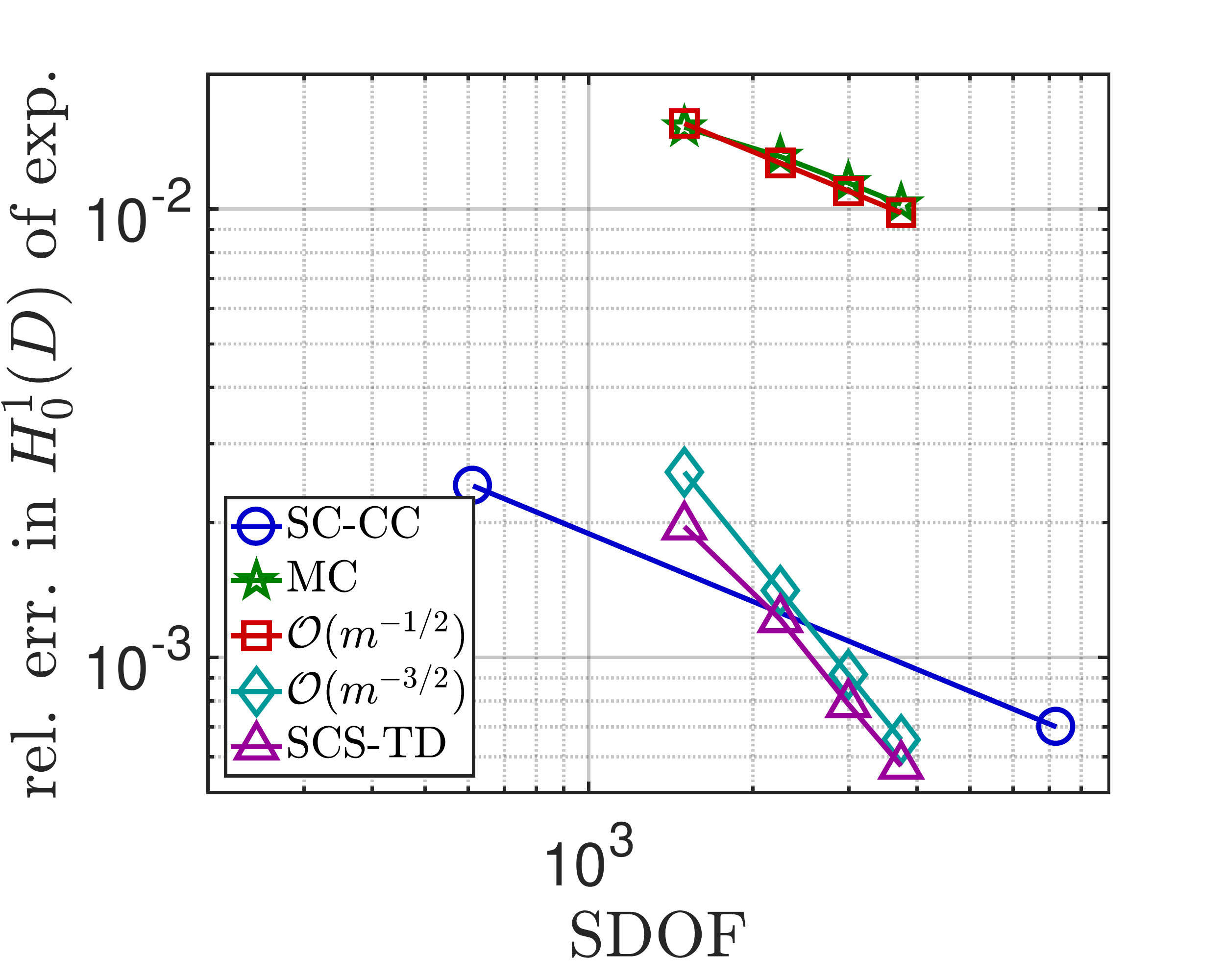}
\includegraphics[clip=true,trim=00mm 04mm 14mm 00mm,width=0.24\textwidth]{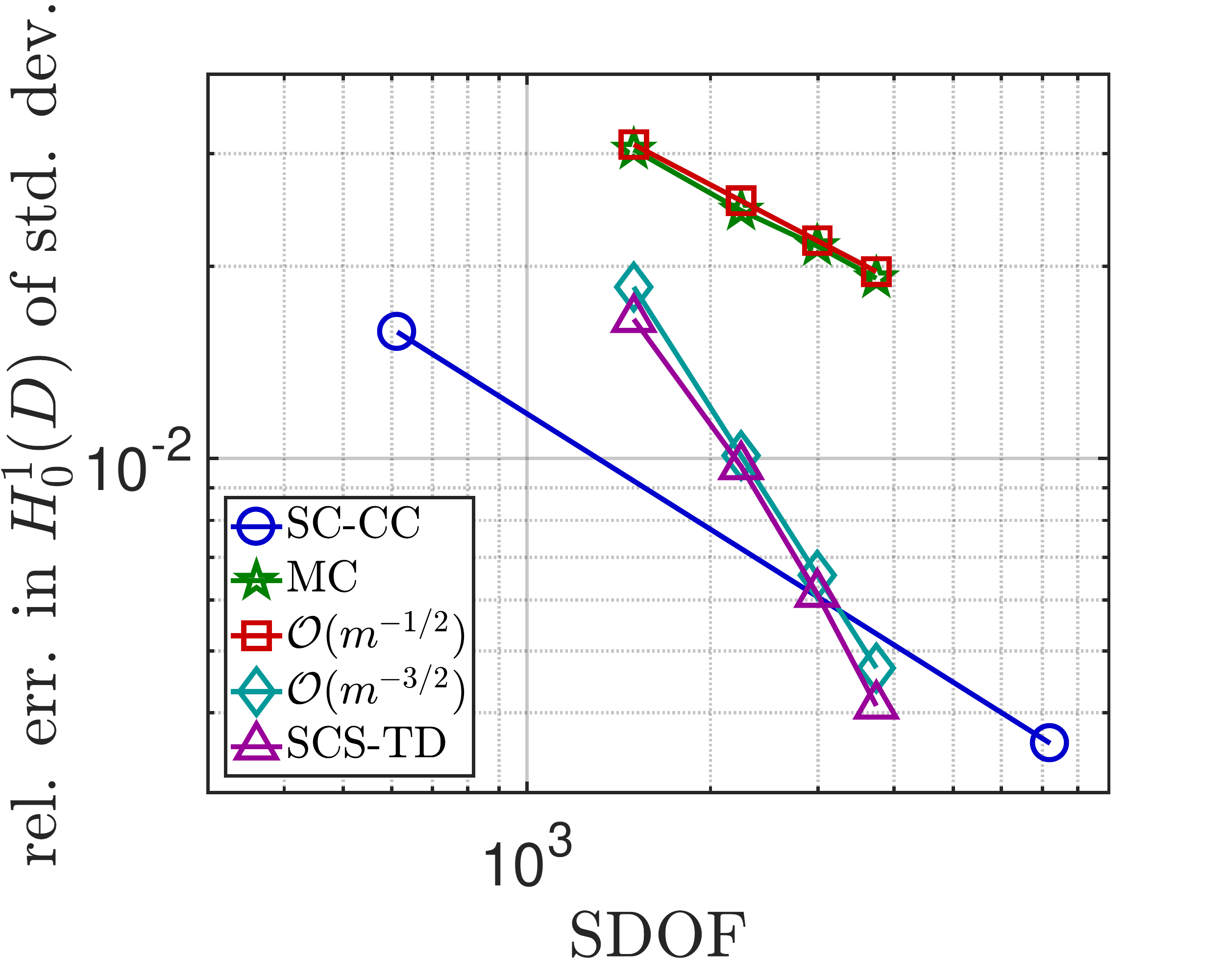}
\end{center}
\caption{Comparison of relative errors in expectation \textbf{(left)} and standard deviation \textbf{(right)} in the $H_0^1(D)$ norm for the $L=2,3$ stochastic collocation (SC-CC), Monte Carlo (MC), and total degree order $p=4$ with $N=5985$ simultaneous compressed sensing (SCS-TD) methods for solving \eqref{eq:model_problem} with coefficient $\log(a(x,\by)-0.5)$ with $a(x,\by)$ from \eqref{eq:transcendental_a_linear_version} and correlation length $L_c=1/8$ in $d= 17$ dimensions. All approximations computed on a finite element mesh with 713 degrees of freedom corresponding to a maximum mesh size of $h\approx 1/32$.}
\label{fig:CS_vs_others_nonlinear_parameterization}
\end{figure}

Figures \ref{fig:CS_vs_others_increasing_correlation_length} \& \ref{fig:CS_vs_others_nonlinear_parameterization} display the achieved results. For highly anisotropic problems, the compressed-sensing based approach is able to naturally detect the underlying anisotropy in refinement. We note that, for the problem under consideration in Figure \ref{fig:CS_vs_others_nonlinear_parameterization}, we are unable to obtain an SG-TD approximation due to the difficulty of solving the nonlinearly coupled SG systems, see \cite{Dexter2016} for more details. We also expect that incorporation of anisotropic weighting schemes such as those considered in \cite{BNTT_comp} and weighted $\ell_1$-minimization techniques such as those from \cite{ChkifaDexterTranWebster18} to improve results for the problems considered. We leave a further study of such improvements to a future work.

\section{Concluding remarks}
\label{sec:conclusion}

In this work we presented {an overview of} a novel sparse polynomial approximation technique, enabling global recovery of solutions to parameterized PDEs.
Our approach, based on extensions of compressed sensing and joint-sparse recovery, treats the solution vector as an element of a tensor product of real Hilbert spaces $\cV$.
The key difference in our technique is the use of a mixed norm involving the energy norm of the associated PDE problem and the standard vector {$\ell_1$} norm.
Within this framework, we are able to prove uniform recovery results through straightforward extensions of concepts such as the restricted isometry and null space properties.
Moreover, by combining extensions of error estimates for the standard basis pursuit denoising problem with quasi-optimal error estimates for sparse approximation of parameterized PDE systems, we are able to derive sub-exponential convergence of our method. For more details see \cite{DexterTranWebster17}.

We have also presented a modification of the standard forward-backward splitting algorithm {for Hilbert-valued recovery}. 
As the considered convex optimization problem is posed over an infinite dimensional Hilbert space, new techniques were needed to establish its convergence properties.
By deriving a novel angular convergence result from the firm nonexpansiveness of the soft-thresholding operator, we are able to prove the strong convergence of the algorithm in the considered setting. 

Finally, we presented numerical results on the application of our approach to the approximation of solutions to both affinely and non-affinely parameterized elliptic PDEs.
We compare our results with those obtained using the stochastic Galerkin, stochastic collocation, and Monte Carlo methods.
The achieved results are positive, highlighting a key benefit of compressed sensing-based approaches, namely the ability to detect underlying problem anisotropy in refinement.

\bibliographystyle{plain}      
\bibliography{database3}   
 
\end{document}